\documentstyle[12pt]{article}
\renewcommand{\baselinestretch}{1.3}

\newtheorem {th}{Theorem}[section]
\newtheorem {lem}[th]{Lemma}
\newtheorem {pr}[th]{Proposition}
\newtheorem {cor}[th]{Corollary}
\newtheorem{defn}{Definition}
\newtheorem{conj}{Conjecture}

\def\Cox{\hfill \Box}
\def\deq{\, {\stackrel {def} {=}}}

\def\dd{\delta}

\def\ee{\epsilon}
\def\Real{\hbox{I\kern-.2em\hbox{R}}}
\def\D{{\bf D}}
\def\DT{\D_\pp (\ppi - I)}
\def\vv{{\bf v}}
\def\VV{{\bf V}}
\def\UU{{\bf U}}

\def\pp{{\bf p}}

\def\SS{{\bf S}}

\def\FF{{\bf F}}
\def\ppi{{\bf \pi}}

\def\ww{{\bf w}}

\def\c1{c_1}
\def\d1{d_1}

\def\sf{\sigma \mbox{-field}}

\def\VT{{\bf{\tilde v}}}
\def\tilt{{\tilde \tau_{n+1}}}
\def\v0{{\bf 0}}

\def\TD{W}
\def\E{{\bf{E}}}
\def\P{{\bf{P}}}

\def\Eb{{\overline {\bf E}}}
\def\Pb{{\overline {\bf P}}}
\def\R{{\bf{R}}}

\def\A{{\cal{A}}}
\def\S{{\cal{S}}}

\def\C{{\cal{C}}}
\def\F{{\cal{F}}}
\def\G{{\cal{G}}}
\def\Fn{{\cal{F}}_n}

\def\B{{\cal{B}}}

\def\nbd{{\cal{N}}}

\def\c1{c_1}
\def\d1{\delta_1}

\def\E{{\bf{E}}}
\def\mv{M_{\vv}\,}
\def\pv{{\pi} ({\vv})}
\def\simp{\bigtriangleup}

\def\TTT{\tilde T}
\def\|{\, | \, }
\def\one{{\bf 1}}

\def\var{\mbox{Var}}

\begin{document}

\begin{titlepage}
\begin{center}
{\large \bf VERTEX-REINFORCED RANDOM WALK} \\
\end{center}

\vspace{5ex}
\begin{quote} \raggedleft
Robin Pemantle \footnote{This research was supported by an NSF graduate
  fellowship and by an NSF postdoctoral fellowship.}\\
Dept. of Statistics \\
U.C. Berkeley \footnote{Now in the department of Mathematics
at the University of Wisconsin-Madison}
\end{quote}

\vfill

{\bf ABSTRACT:} \break

\renewcommand{\baselinestretch}{1.0}\large\normalsize
\noindent{This paper considers} a class of non-Markovian discrete-time 
random processes on a finite state space $\{ 1 , \ldots , d \}$.  The 
transition probabilities at each time are influenced by the number of times
each state has been visited and by a fixed a priori likelihood matrix, $\R$,
which is real, symmetric and nonnegative.
Let $S_i (n)$ keep track of the number of visits to state $i$ up to time $n$,
and form the fractional occupation vector, $\VV (n)$, where 
$v_i (n) = S_i (n) / (\sum_{j=1}^d S_j (n) )$.  It is shown
that $\VV (n)$ converges to to a set of critical points for the
quadratic form $H$ with matrix $\R$, and that under nondegeneracy conditions
on $\R$, there is a finite set of points such that with probability one,
$\VV (n) \rightarrow \pp$ for some $\pp$ in the set.
There may be more than one $\pp$ in this set for which $\P (\VV (n)
\rightarrow \pp) > 0$.  On the other hand $\P (\VV (n) \rightarrow \pp)
= 0$ whenever $\pp$ fails in a strong enough sense to be maximum for $H$. 
\renewcommand{\baselinestretch}{1.5}\large\normalsize

\end{titlepage}

\section{Introduction}

This paper considers a stochastic process in discrete time on a finite
state space $\{ 1 , \ldots , d \}$, in which the probability of a 
transition to site $j$ increases each time $j$ is visited.  To define
the process, let $\R$ be a real symmetric $d \times d$ matrix with $\R_{ij}
\geq 0$ for each $i,j$, and $\sum_i \R_{ij} > 0$ for each $j$.
For $n \geq d$, inductively define random variables $Y_n$ and
$\SS (n) = (S_1 (n) , \ldots , S_d (n))$ as follows.  Let $S_i (d) = 1$
for $i = 1 , \ldots , d$ and let $Y_d = 1$.  Let $\F_n$ be the $\sf$
generated by $Y_j : d \leq j \leq n$ and let $Y_{n+1}$ satisfy
$$ \P (Y_{n+1} = j \| \F_n) = \R_{Y_n , j} S_j (n)
    / \sum_i \R_{Y_n , i} S_i (n) .$$
Let $S_i (n+1) = S_i (n) + \dd_{Y_{n+1},i}$.  In other words,
$\SS (n)$ counts one plus the number of times $Y$ has occupied each 
state.  The sequence of ordered pairs $(Y_n , \SS (n))$ is a 
Markov chain, whereas the sequence $Y_n$ is not.  

Define $\VV (n) = \SS (n) / n$, so that each $\VV (n)$ is an element
of the $d-1$-simplex $\simp \subseteq \Real^d$.  (In general,
boldface is used for vectors and lightface is used for their components.)  
This paper studies
the question of when $\VV (n)$ converges and to which possible limits.
Since $\VV (n)$ may be viewed as an empirical occupation measure for the 
$Y$ process, this is essentially asking whether $Y$ obeys a strong
law of large numbers.  A few remarks about the model are in order.

The process is meant to model learning behavior.  Think of $\R_{ij}$
as a set of initial transition probabilities; each time $Y$ visits
site $j$, this choice is positively reinforced, resulting in
transition probabilities proportional to $\R_{ij} S_j$.  The choice
of starting state, $Y_d = 1$, is arbitrary; also, setting
each $S_i (d)$ equal to one is a matter of convenience and 
in fact the theorems in this paper are true for any choice of
$S_i (d) > 0$ and any $Y_d \in \{1 , \ldots , d\}$.  The 
requirement that $\R$ be symmetric may not always be reasonable
in applications, but is essential for our arguments.

Similar models have been studied in \cite{IT} under the name of {\em 
random processes with complete connections}.  When the entries of $\R$ 
are all one, the model reduces to a P\'olya urn model; the behavior
in this case is atypical, since most of our results apply to the
``generic'' case where $\R$ is invertible.  Another similar
process called {\em edge-reinforced random walk} is studied in
\cite{Di,Pe1,Pe2,Da}; in that case, transitions from $i$ to $j$
are positively reinforced each time a transition is made from
$i$ to $j$ or $j$ to $i$.  Thinking of the process as traversing
a graph with vertices $1 , \ldots , d$, this kind of reinforcement
keeps track of moves along each edge of the graph, while the process 
studied in the present paper keeps track of visits to each vertex.
Strong laws for edge-reinforced random walk can be found in 
\cite{Di,Pe1,Da}.  

The remainder of this introductory section motivates and states
the main results.  Subsequent sections give proofs of of the four
results.  Examples and open questions are discussed in the final section.
\begin{defn} 
For $\vv \in \simp$, let $N_i (\vv ) = \sum_j \R_{ij} v_i$.  Abbreviate
this by $N_i$ when a particular vector $\vv$ may be understood.
\end{defn}
\begin{defn}
For $\vv \in \simp$, let $H(\vv ) = \sum_i v_i N_i (\vv ) = \sum_{ij}
\R_{ij} v_i v_j$.
\end{defn}
\begin{defn} \label{eqpi}
For $\vv \in \simp$ such that $H(\vv ) > 0$, define a vector
$\ppi (\vv ) \in \simp$ by $\pi_i (\vv ) = v_i N_i (\vv ) / H( \vv)$.
\end{defn}
\begin{defn}
For $\vv \in \simp$ such that $H(\vv ) > 0$, define a Markov transition
matrix $M(\vv )$ by $M_{ij} (\vv ) = \R_{ij} v_j / N_i$.  
\end{defn}
Note that $H(\VV (n))$ is 
below by $\min \{ \R_{ij} : \R_{ij} > 0 \}$.  Thus $H$ never
vanishes on the closure of the set of possible values of $\VV (n)$,
and the clauses about $H$ not vanishing in the above definitions
are merely {\em pro forma}.  For a fixed $\vv$,
$(\ppi M)_i = \sum_j \pi_i M_{ij} = \sum_j (v_i N_i/H) (\R_{ij} 
v_j/N_i) = \sum_j v_i v_j \R_{ij} / H = v_i N_i /H = \ppi_i$, so
$\ppi (\vv )$ is an invariant probability for the transition matrix
$M(\vv)$.  The behavior of $\VV (n)$ can heuristically be explained
as follows.

For $n \gg L \gg 1$, compare $\VV (n+L)$ to $\VV (n)$.  Since $n \gg L$,
the $Y$ process between these times behaves as if $\VV$ is not changing,
and hence approximates a Markov chain with transition matrix 
$M(\VV (n))$.  Since $L \gg 1$, the occupation measure between
these times will be close to the invariant measure $\ppi (\VV(n))$.
This means that $\VV(n+L) \approx \VV(n) + (L/n) (\ppi (\VV(n)) - \VV(n))$.
Passing to a continuous time limit gives
\begin{equation} \label{flow}
{d \over dt} \, \VV (t) = {1 \over t} (\ppi (\VV (t)) - \VV (t)) .
\end{equation}
Up to an exponential time change, $\VV$ should then behave like an
integral curve for the vector field $\ppi - I$.  One would expect
convergence to a critical point or set and, because of the random 
perturbations, one would not expect convergence to any unstable
equilibrium.  It is not in general possible to find a potential for this
vector field, but the function $H$ is a Lyapunov function for it.
Then one expects convergence of $\VV (n)$ to a maximum for $H$.  
\begin{defn}
Let $\C \subseteq \simp$ be the set of points $\vv$ for which $\ppi 
(\vv) = \vv$.  The term {\em critical point} will be used to denote
points of $\C$.  Let $\C_0 \subseteq \simp$ bet the set of points $\vv$
for which $M(\vv)$ is reducible.
\end{defn}

Section~\ref{preliminaries} will discus the nature of $\C$ and $\C_0$, 
and give conditions under which Theorem~\ref{th converge} 
(proved in Section~\ref{pf converge}) implies 
almost sure convergence of $\VV (n)$. 
\begin{th} \label{th converge}
With probability one, $dist(\VV (n) , \C \cup \C_0) \rightarrow 0$, where 
$dist(x,A)$ denotes $\inf \{ |x - y| : y \in A \}$.
\end{th}
\begin{defn}
For $\vv \in \simp$, define $face (\vv) = \{ \ww \in \simp : \forall i , 
v_i = 0 \mbox{ implies } w_i = 0 \}$ to be the closed face of
$\simp$ to which $\vv$ is interior.
\end{defn}
\begin{defn}
For any $\pp \in \C$ that is in a proper face of $\simp$ a {\em linear
non-maximum} iff 
\begin{equation} \label{lnm}
\D_\pp H(e_k - e_j) > 0 \mbox{ for some } e_k \notin
face(\pp) , e_j \in face (\pp) .
\end{equation}
(Here $e_1 , \ldots , e_d$ are the standard basis vectors in $\Real^d$.)
\end{defn}
The following theorems, proved in Section~\ref{pf nonconverge2}
and~\ref{pf nonconverge1} respectively, give conditions under which 
convergence to a critical point is impossible.
\begin{th} \label{interior nonconvergence}
Suppose that $\R$ is nonsingular and let $\pp$ be the unique critical
point in the interior of $\simp$.  Then $\P (\VV (n) \rightarrow \pp)
= 0$ whenever $\pp$ fails to be a maximum for $H$.  This happens
if and only if $\R$ has more than one positive eigenvalue, which happens
if and only if the linear operator $\DT$ on $-\pp + \simp$ has a 
positive eigenvalue.
\end{th}
\begin{th} \label{bdry nonconvergence}
Suppose $\pp$ is a linear non-maximum in a proper face of $\simp$.  
Then $\P (\VV (n) \rightarrow \pp) = 0$. 
\end{th}
A sort of converse to these nonconvergence theorems gives a criterion
for convergence with positive probability of $\VV (n)$ to stable
critical points.  This is proved in Section~\ref{pf converge}
\begin{th}  \label{th positive prob converge}
Let $A$ be a component of $\C$ disjoint from $\C_0$ and
suppose that $A$ is a local maximum for $H$ in the sense that 
there is some neigborhood $\nbd$ of $A$ for which $\vv \in \nbd$
and $\pp \in A$ imply $H(\vv ) < H(\pp)$.  Then $\P(dist(\VV (n) , A)
\rightarrow 0) > 0$.
\end{th}

\section{Preliminaries} \label{preliminaries}

The following proposition verifies that $H$ is a Lyapunov function
for the vector field $\ppi - I$ and gives alternate characterizations
of the set of critical points.  The notation used throughout
for vector calculus is $\D_\vv F(\ww)$ to denote the derivative
of $F$ in the direction $\ww$ at the point $\vv$, thus $\D_\vv F$
denotes the linear operator approximating $F(\vv + \cdot) - F(\vv)$.

\begin{lem} \label{Hincr}
For any $\vv \in \bigtriangleup \, , \, \D_\vv H(\pv -  \vv ) \geq 0$.
Furthermore, the following are equivalent:
\begin{equation}
\begin{array} {cl}
(i) & \D_\vv H(\pv - \vv ) = 0 \\
(ii) & \D_\vv H |_{\mbox{face} (\vv )} = \vec{0} \\
(iii) & \mbox{for those i such that } v_i > 0 \, , \, N_i \mbox{ are equal}
       \\
(iv) & \mbox{for all } i , \, v_i = \sum_j \R_{ij} v_i v_j / N_j \\
(v)  & \pi (\vv ) = \vv
\end{array} \label{crit}\end{equation}
where $0/0 = 0$ in (iv) by convention. 
\end{lem}

\noindent{Proof:}  For fixed $i$ and $j$ and 
constant $c$, consider the operation of increasing $v_j$ by the 
quantity $c v_i v_j (N_j - N_i)$ and decreasing $v_i$ by the same amount.
When $c = 1 / H(\vv )$ and this operation is done simultaneously for every
(unordered) pair $i$, $j$, then the resulting vector is $\pv$:
the next value of the $i^{th}$ coordinate is given by 
$$\begin{array}{ll} & v_i + (1 / H(\vv )) (\sum_j v_i v_j N_i - 
    \sum_j v_i v_j N_j) \\
= & v_i + (1 / H(\vv )) (v_i N_i - v_i H(\vv )) = \pi_i (\vv ) .  \end{array} $$
So an infinitesimal move towards $\pv$ corresponds to doing these
additions and subtractions simultaneously with an infinitesimal $c$.  To 
show that this increases $H$, it suffices to show that for each unordered
pair $i$, $j$, the value of $H$ is increased, since $H$ is smooth and therefore
well approximated by its linearization near any point.  So let $i$, $j$ be 
arbitrary.  Writing $\vv^{(1)}$ for the new vector gives
\begin{eqnarray*}
H(\vv^{(1)}) & = & \sum \R_{rs} {v_r}^{(1)} {v_s}^{(1)} \\
& = & \sum_{r,s} \R_{rs} v_r v_s + 2 \sum_s \R_{is} cv_i v_j 
(N_i - N_j ) v_s \\
& & + 2\sum_r \R_{rj} c v_i v_j (N_j - N_i ) \\
& = & H(\vv ) + 2 c v_i v_j (N_i - N_j )^2 \\
& \geq & H(\vv ) 
\end{eqnarray*}
so $H$ is nondecreasing.  This proves the first part.  

For the equivalences, first note that
if there are any $i$ and $j$ for which $N_i \neq N_j$ and neither 
$v_i$ nor $v_j$ is zero, then $H$ strictly increases.  
Thus $(i) \Leftrightarrow (iii)$.  Since 
\begin{equation} \label{nablaH} \D_\vv H \mbox{ is just inner product with the
vector }(2N_1, \cdots , 2N_n),   \end{equation}
and restricting to $\mbox{face}
(\vv )$ just throws out the coordinates $i$ such that $v_i = 0$, it is
easy to see that $(ii) \Leftrightarrow (iii)$.  Assuming $(iii)$, suppose
the common value of the $N_i$ is $c$.  Then multiplying $(iv)$ by $c$
gives $\sum_j v_i v_j = c \cdot v_i$, so $(iii) \Rightarrow (iv)$.
Now assume $(iv)$.  Letting $\mv$ denote the matrix as well as the
Markov chain, $(iv)$ just says that $\vv$ is stationary for $\mv$.
Then $\pv - \vv = \vec{0}$ so $(v)$ holds.   And finally, $(v)
\Rightarrow (i)$ trivially.    $\Cox$

\begin{pr} \label{critset}
The set $\C$ has
finitely many connected components, each of which is closed and on
each of which $H$ is constant.  Furthermore, if all the principal 
minors of $\R$ are invertible, then $\C$ consists of at most 
$2^d - 1$ points. 
\end{pr}

\noindent{Proof:}  By (\ref{crit}) (ii), $\C$ is the union over 
all $2^d - 1$ faces
$F$ of the sets $\C_F = \{\vv : \D_\vv H |_F (\vv ) = 0 \}$.  By~(\ref{nablaH})
and the comment following, $\D_\vv H |_F$ is linear, so $\C_F$ is a closed,
convex, connected set.  It is easy to see that $H$ is constant on $\C_F$ by
integrating $\D_\vv H |_F$.  The first part of the proposition follows 
since each connected component of $\C$ is the union of some of the $\C_F$. 
For the second part, fix a face $F$ and let $\R_F$ be the matrix gotten
from $\R$ by deleting rows and columns indexed by those $i$ for 
which $v_i = 0$ for all $\vv \in F$.  If this is invertible, then 
equation~(\ref{crit}) (iii) implies that the only possible element of
$\C$ in the interior of $F$ is whichever multiple of $(1 , \ldots , 1)
\R_F^{-1}$ lies on the unit simplex.   $\Cox$ 

If all the off-diagonal entries of $\R$ are positive, 
it is immediate that $M(\vv)$
is irreducible for all $\vv \in \simp$.  Conversely, if 
$\R_{ij} = 0$ for some $i \neq j$, then $M(\vv)$
is reducible when $\vv$ is any nontrivial combination of
$e_i$ and $e_j$.  Thus it a necessary and sufficient condition for
$\C_0$ to be empty is that $\R_{ij} > 0$ off of the diagonal.
In any event, $\C_0$
is a union of proper faces of $\simp$.  The following corollary
to Theorem~\ref{th converge} is now immediate.

\begin{cor} \label{converge to point}
If all the off-diagonal entries of $\R$ are positive and all the principal 
minors of $\R$ are invertible, then $\VV(n)$  converges almost surely.
\end{cor}
$\Cox$

\section{Proofs of convergence results} \label{pf converge}

The proof of Theorem~\ref{th converge} begins with a lemma giving
a lower bound on the expected growth of $H(\VV(n))$ when 
$\VV(n)$ is not near $\C \cup \C_0$.  

\begin{lem} \label{lem4.1}
Let $\nbd$ be a closed subset of the simplex, with $\nbd \cap (\C \cup 
\C_0 ) = \emptyset$.  Then there exist an $N$, $L$ and $c > 0$
such that for any $n > N$, $\E (H(\VV (n + L)) \| \VV (n)) > H(\VV (n)) 
+ c / n$ whenever $\VV (n) \in \nbd$.  
\label{submart} \end{lem}

\noindent{Proof:}  For any $n$, let $M_n (n), M_n (n+1) , \ldots$ denote 
a Markov chain beginning at $Y_n$ at time $n$, whose transition
matrix thereafter does not change with time and is given by $M(\VV (n))$.  
Let $\SS' (n) = \SS (n)$ and for $i > n$, let $\SS' (i) = \SS' (i-1)
+ e_{M_n (i)}$, where $e_j$ is the $j^{th}$ standard basis vector.
Let $\VV' (i) = \SS' (i) / i$.

First I claim that the lemma is true with the Markov process $\VV'$ 
substituted for $\VV$.  
By Lemma~\ref{Hincr}, $\D_\vv H(\pi (\vv ) - \vv ))$ is nonzero on
$\nbd$, so by compactness it is bounded below by some $c_0$ on $\nbd$.  
Choose any $c_1 < c_0$.  The occupation measure of a process between 
times $N$ and $N+L$ can change by at most $L/(N+L)$ in total variation.
Since $H$ is smooth, it is possible to choose
$N/L$ large enough so that whenever $n \geq N$,
$H[\VV ' (n) + (L/(N+L))(\ppi (\VV  (n)) - \VV (n))] > c_1 L / (n+L)$. 
By the Markov property, $(\SS' (n + L) - \SS (n)) / L$ approaches a point-mass
at $\ppi (\vv )$ in distribution as $L$ increases.  In fact, the rate of
convergence of $M^k (\VV(n )) \ww$ to $\ppi (\VV(n ))$ 
is exponential and controlled by the second-largest eigenvalue 
of $M(\VV(n))$ according to the Perron-Frobenius theorem.  If $M(\VV (n))$
is aperiodic, then since $M(\vv )$
varies continuously with $\vv$, eigenvalues
are continuous, and the non-degeneracy hypothesis says that $\nbd$
contains no points where the second-largest
eigenvalue is 1, the second-largest eigenvalue is bounded away from 1.  
It follows that a large enough $L$ may be chosen uniformly in $\vv$
so that $\E (H(\VV' (n + L)) - H(\VV' (n)) \| \F_n) > c / n$ for any 
$c_2 < c_1$, and the claim is established.  If $M(\VV (n))$ is periodic, 
then it has period 2 and a simple eigenvalue at $-1$; the claim
follows in this case from grouping together pairs of times
$2n$ and $2n+1$.

Now couple the Markov chain $\VV' (n+i)$ to $\VV (n+i)$ in such a 
way so the two 
move identically for as long as possible.  Formally, define $\{M_n (i) \}$
and $\{ Y_i \}$ on a common measure space so that if $Y_j = M_n (j)$
for all $n < j < n+k$ then 
$$\P (Y_{n+k} \neq M_n (n+k) ) \| Y_{n+k-1} = i) = \sum_j
{1 \over 2}\, |M_{ij} (\VV(n+k)) - M_{ij} (\VV (n))| . $$
Picking $c < c_2$ and $N/L$ large enough so that \begin{equation}
(L^2 / N) (L/N) || \D H ||_{op} < (c_2 - c) / N ,
\label{eq1} \end{equation} 
the coordinates of $\VV$ 
cannot change by more than $L/N$ in $L$ steps, so the probability of
an uncoupling at any of the $L$ steps is bounded by $L^2 / N$. 
Then $\E |H(\VV (n+L)) - H(\VV' (n+L))|
< (c_2 - c)/N$ by~(\ref{eq1}), and combining this with the earlier claim
proves the lemma.    $\Cox$

Before proving Theorem~\ref{th converge}, here is a sketch of the argument.
On any set $\nbd$ away from $\C \cup \C_0$, Lemma~\ref{submart} says
the expected value of $H(\vv (n))$ grows, provided you sample at time intervals
of size $L$.  The cumulative differences
between $H(\vv (n+L))$ and $\E (H(\vv (n+L)) \| \vv (n))$ form a convergent
martingale, so $H(\VV (n))$ itself is growing at rate $c/n$ when $\VV (n)
\in \nbd$.  The rate of change in position of $\VV (n)$ is also order $1/n$
per step, so if $\VV$ goes from one given point
of $\nbd$ to another, $H(\VV (n))$ increases by an amount independent of 
time.  The only way it can decrease again is for $\VV (n)$ to leave $\nbd$
at a place where $H$ is large and re-enter where $H$ is small.  The effect of
such a possibility can be made arbitrarily small because $H$ is 
nearly constant on the connected components of $\simp \setminus \nbd$.

\noindent{Proof:} of Theorem~\ref{th converge}:  Since
the connected components, $\C_i ,\ldots \C_k$ of 
$\C \cup \C_0$ are closed,
$m = \min \{ d(\C_i , \C_j ) \} > 0$.  Pick any $r < m/3$.  Let
\begin{eqnarray} 
{\nbd_1}^i & = & \{ \vv : d(\vv , \C_i ) < r \} \label{n1i} \nonumber \\
\nbd_1 & = & \simp \setminus\bigcup_{i=1}^k {\nbd_1}^i . \label{n1}
\end{eqnarray}
Note that 
\begin{equation} \label{fatenough}
i \neq j \Rightarrow d({\nbd_1}^i , {\nbd_1}^j ) > r .
\end{equation}
By the preceding lemma with $\nbd = \nbd_1$, $c_1 , L_1 , N_1$ can be found
for which $n \geq N_1 $ implies
$\E (H(\VV (n+L)) \| \VV (n)) \geq H(\VV (n)) + c/n$.  Pick any $L' > L_1 $ and
define
\begin{eqnarray*} 
{\nbd_2}^i & = & {\nbd_1}^i \cap  \{ \vv : |H(\vv ) - H( \C_i )| 
   < rc/2L' \} %%%\label{n2i} 
   \\
\nbd_2 & = & \simp \setminus\bigcup_{i=1}^k {\nbd_2}^i . %%%\label{n2}
\end{eqnarray*}
Figure 1 gives an example of these definitions when $d = 3$; the
heavy lines are the boundary of $\nbd_1$ and the lighter lines are the boundary
of $\nbd_2$.  

Apply the lemma to $\nbd_2$ to get $N_2 , c_2$ and $L_2$.  Define the process
$\{ \UU (n) \}$ that samples $\VV (n)$ at intervals of $L_1$ 
on $\nbd_1$ and $L_2$ elsewhere, by
\begin{eqnarray*}
\UU (n , \omega ) & = & \VV (f(n , \omega )) \\
& & \\ & \mbox{ where } & \\ & & \\
f(1 , \omega ) & = & \max \{ N_1 , N_2 \} \mbox{	and } \\
f(n+1 , \omega ) & = & 
  \left \{ \begin{array}{l} f(n , \omega ) + L_1 
    \mbox{ if } \VV (f(n , \omega )) \in \nbd_1 ; \\ 
  f(n , \omega ) + L_2 \mbox{ if } 
    \VV (f(n , \omega )) \notin \nbd_1 . \end{array} \right. .
\end{eqnarray*}
Clearly, $\UU (n)$ converges if and only if $\VV (n)$ converges.  Letting
$U(n) = H(\UU (n))$, write $U(n) = M(n) + A(n)$ where $\{M(n)\}$
is a martingale and $\{A(n)\}$ is a predictable process with respect to 
$\F_{f(n)}$.  The key properties needed are
\begin{eqnarray}
& & M(n) \mbox{ converges almost surely} \label{mconv} \\
& & A(n+1) \geq A(n) + c/n \mbox{ if } \UU (n) \in \nbd_1 \label{agrow} \\
& & A(n+1) \geq A(n) \mbox{ if } \UU (n) \in \nbd_2 . \label{aweakgrow} 
\end{eqnarray}
To verify (\ref{mconv}), note that $|U(n+1) - U(n)| \leq 
\max \{ L_1 , L_2 \} / f(n) = O(1/n)$, since by~(\ref{nablaH}), $H$ is 
Lipschitz on $\simp$.  Then $|M(n+1) - M(n)| =
O(1/n)$ as well, so $M(n)$ converges in $L^2$, hence almost surely.
Properties~(\ref{agrow}) and~(\ref{aweakgrow}) are evident from 
the construction.

The next thing to show is Claim 1: $\UU (n) \in {\nbd_2}^a$ infinitely 
often for at most one $a$ almost surely.  Consider any sample path $\UU
(1) , \UU (2) ,
\ldots$.  For $n < t$, define the event $\B (a,b,n,t,\omega )$ to occur if
\begin{equation} \label{bad}
\UU (n) \in {\nbd_2}^a \mbox{ and } \UU (t) \in {\nbd_2}^b \mbox{ with }
\UU (i) \in {\nbd_2} \mbox{ for all $i$ such that } n<i<t .
\end{equation}
If $\B (a,b,n,t,\omega )$ occurs, let 
\begin{eqnarray*}
r & = & \max \{ i : n \leq i < t \mbox{ and } \UU ( i) \in {\nbd_1}^a \} \mbox{ and } \\[2ex]
s & = & \min \{ i : n \leq i < t \mbox{ and } \UU ( i) \in {\nbd_1}^b \}
\end{eqnarray*}
be respectively the last exit time of ${\nbd_1}^a$ and the first entrance time
of ${\nbd_1}^b$.
The dotted path in figure 1 gives an example of this.  
By~(\ref{agrow})
and~(\ref{aweakgrow}),
$$ A(i+1) - A(i) \geq c/i \mbox{ for } r < i < s$$
$$ A(i+1) - A(i) \geq 0 \mbox{ for } n < i < t .$$
Then 
\begin{eqnarray*}
& & A(t) - A(n) \\
& = & [A(t) - A(s)] + [ A(s) - A(r+1)] \\
& & + [A(r+1) - A(n+1)] + [A(n+1) - A(n)] \\
& \geq & 0 + \left ( \sum_{i=r+1}^{s-1} c / i \right )+ 0 - L_2 / n \\
& = & O(1/n) +  (c/L_1 ) \sum_{i = r}^{s-1} L_1 / i \\
& \geq & O(1/n) + (c/L_1 ) \sum_{i = r}^{s-1} |\UU (i+1) - \UU (i)| \\
& \geq & O(1/n) + (c/L_1 ) |\UU (s) - \UU (r)| \\
& > & O(1/n) + rc/L_1 
\end{eqnarray*}
by (\ref{fatenough}).  Now $U(t) - U(n) \leq H(\C_b ) - H(\C_a ) + rc/L'$
by the construction of $\nbd_2$.  So $M(t) - M(n) \leq H(\C_b ) - H(\C_a )
+ rc/L' - rc/L_1 + O(1/n)$.  If $H(\C_b ) \leq H(\C_a )$, the choice of $r$
guarantees that this expression
is strictly negative and bounded away from 0 for large $n$.  Therefore
if $M(n)(\omega )$ converges, then $\B (a,b,n,t,\omega )$ happens only
finitely often for $a,b$ such that $H(\C_b ) \leq H(\C_a )$.  But then 
it happens only finitely often for any $a \neq b$, since $\UU$ can make
only $k-1$ successive transitions from ${\nbd_2}^a$ to ${\nbd_2}^b$ with
$H(\C_b ) > H(\C_a )$.  Thus the almost sure convergence of $M(n)$ implies
that $\UU (n) \in {\nbd_2}^a$ infinitely often for at most one $a$
almost surely and Claim 1 is shown.

In other words, transitions between small neighborhoods of $\C_i$ and $\C_j$
eventually cease for $i \neq j$.  Claim 2 is that $\VV (n)$ may not oscillate
between a small neighborhood of $\C_i$ and a set bounded away from $\C$.
To show this, require now that
$r < m/6$.  With $\nbd_1$ and $\nbd_2$ defined as before, 
define $\nbd_3 \subseteq \nbd_1$ by (\ref{n1}) with $2r$ in place of $r$.  
Since $2r < m/3$, 
equation~(\ref{fatenough}) holds with $\nbd_3$ in place of $\nbd_1$.
An argument identical to the one establishing Claim 1 now shows that with 
probability 1 there
are only finitely many values of $n$ and $t$ for which
$$\UU (n) \in {\nbd_2}^a , 
\UU (i) \in {\nbd_3} \mbox{ and } \UU (t) \in {\nbd_2}^a \mbox{ for } n<i<t.$$
[The argument again: $A(i)$ is nondecreasing when $\UU (n) \in {\nbd_1}^a$ 
and increases by at least the fixed amount $rc/L_1$ each time $\UU$ makes the 
transit from ${\nbd_1}^a$ to $\nbd_3$.  The increase in $A$ is greater than
the greatest difference in values of $H$ taken at two points of ${\nbd_2}^a$,
so the martingale $M$ must change by at least $rc/L_1 - rc/L'$ during every
transit.  Since $M$ converges, this happens finitely often.] 

Claim 3 is that the event $\{ \omega : \UU ( t,\omega ) \in \nbd_1 
\mbox{ for all }t>n \}$ has probability 0 for each $n$; it is proved
in an identical manner.  Putting together Claims 1 and 3, it follows that
for any small $r$ there is precisely one $a$ for which $\UU (n) \in 
{\nbd_1}^a$ infinitely often.  Then by Claim 2 for a different $r$, 
$\nbd_3$ stops being visited, so 
letting $r \rightarrow 0$ proves the theorem.   $\Cox$

The proof of Theorem~\ref{th positive prob converge}
is just an easier version of the proof of Theorem~\ref{th converge}.  

\noindent{Sketch} of proof of Theorem~\ref{th positive prob converge}:
A process $\UU(n)$ may be defined as in the previous proof, so 
that $\VV (n)$ converges iff $\UU(n)$ converges and so that
$U(n) \deq H(\UU (N))$ breaks into a martingale $M(n)$ and
a predictable process $A(n)$.  Note that the argument showing
an $L^2$ bound of $c/n$ on $M(\infty) - M(n)$ still works
conditionally on $\UU (n)$.  By a standard maximal inequality,
given any $\ee > 0$, an $n$ may be
chosen large enough so that $\P ( \inf \{M(n) - M(n+i) : i > 0 \} < -\ee \|
\UU (n)) < \ee$.  The assumptions of the theorem imply the
existence of an $\ee$ for which the component $B$ of $H^{-1} [a-2\ee ,
a]$ is disjoint from $(\C \cup \C_0) \setminus A$, where $a$ is
the value of $H$ on $A$.  Now for sufficiently large $n$, the
event $U (n) \in H^{-1} [a-\ee,a] \cap B$ has positive probability.
Conditional on this event, the probability that $M(n+i) -M(n)$
never goes below $-\ee$ has been shown to be less than $\ee$
for large $n$.  Since $dist(U (n) , \C \cup \C_0) \rightarrow 0$
by Theorem~\ref{th converge}, and $\UU (n)$ cannot leave $B$
without $U(n)$ becoming less than $a - 2\ee$, it follows that
$dist(\UU (n) , A) \rightarrow 0$, proving the theorem.    $\Cox$

\section{Proof of Theorem~{\protect \ref{bdry nonconvergence}}} 
\label{pf nonconverge1}

To prove Theorem~\ref{bdry nonconvergence}, begin by seeing why
it should be true.  With $\pp$ as in the statement of the theorem,
equation~(\ref{crit})~(iii) says that
the $N_i$ have a common value, $\lambda$, for those $i$ such that $p_i > 0$. 
Assuming~(\ref{lnm}) for a given $e_k$ and using 
equation~(\ref{nablaH}) for $\D H$ shows that $N_k > N_j = \lambda$.  So
\begin{equation} \label{ineq}
\sum_i \R_{ki} p_i / N_i = \sum_{p_i > 0} \R_{ki} p_i / \lambda = N_k / 
\lambda = 1 + b
\end{equation}
for some $b > 0$, $k$ such that $p_k = 0$.  Now when $\VV (n)$ is close
to $\pp$, $v_k (n)$ will be close to but not equal to zero.  The
expected number of visits to state $k$ during a period of time 
from $n$ to $n+T$ in which the occupation measure is close to $\pp$ 
will be approximately $T \sum_i p_i (\R_{ik} v_k / N_i) = 
T v_k N_k / \lambda = (1+b) T v_k$.  In other words, $v_k$ 
will begin to increase and $\pp$ should be an unstable point with 
no possibility of $\VV (n)$ converging there.  The actual proof will
consist of making this rigorous.
 
To avoid bogging down in trivialities, $\SS (n)$ and $\VV(n)$ will be used
to stand for $\SS (\lfloor n \rfloor )$ and $\VV (\lfloor n \rfloor )$.
Inequalities will be verified as if $n$ were an integer;
it is always possible to choose epsilons and deltas a little bit
smaller to compensate for the roundoff errors.  Begin by recording 
a few propositions whose proofs are omitted when elementary.

\begin{pr} \label{convex}
Fix $\pp$ and let $\nbd_1$ be a neighborhood of $\pp$.  For any $\dd > 0$ there
is a neighborhood $\nbd$ of $\pp$ included in $\nbd_1$ such that for all 
$n > 1 / \dd$, the two conditions
\begin{eqnarray*}
(i) & \VV (n) \in  \nbd & \mbox{ and} \\ 
(ii)& \VV (n+\dd n) \in \nbd & 
\end{eqnarray*}
imply
\begin{eqnarray*} (iii) & (\SS (n+\dd n ) - \SS (n))/\dd n \in \nbd_1 & . 
\end{eqnarray*}
$\Cox$
\end{pr}

The heuristic calculation at the beginning of this section
is made precise as follows.
\begin{pr} \label{1+b}
Let $\pp ,k,b$ be such that~(\ref{ineq}) holds and let $\SS$ be any
vector function of $n$.  Then there is an $\ee > 0$
and a neighborhood $\nbd_1 = \{ \vv \in \simp : |\vv - \pp | < \ee \}$ such
that for all $\dd > 0$ and for all $n$, the conditions $\VV (n) \in \nbd_1$
and $(S_i (n + \dd n) - S_i (n))/ \dd n \geq p_i - \ee$ for all $i$ imply
\begin{equation} \label {suflarge}
\sum_i (S_i (n + \dd n) - S_i ( n))\R_{ik} v_k (n) / (1 + \dd )N_i (n)
> \dd {1 + b/2 \over 1+\dd} S_k (n) .\end{equation}
\end{pr}

\noindent{Proof:}  As $\ee \rightarrow 0, \,1/n$ times the left-hand 
side converges to
$\dd p_k (n) / (1 + \dd ) \sum_i p_i (n) \R_{ik} / N_i (n)$ $= \dd p_k (n)
(1 + b) / (1 + \dd )$ while $1/n$ times the right-hand side converges to 
$\dd p_k (n) (1 + b/2) / (1 + \dd )$.  Since the convergence is uniform in 
$\dd$, the result follows.   $\Cox$

\begin{pr}
\label{pointmass}
Let $b > 0$ and $\ee_1 > 0$ be given.  Let $\{ B_{\alpha} \}$ be a collection
of independent Bernoulli random variables with $\E ( \sum_{\alpha} B_{\alpha} )
\geq (1+b)L$.
There exists an $L_0$ such that whenever $L > L_0$, 
$\P (\sum_{\alpha} B_{\alpha} / L > 1 + b/2 ) > 1 - \ee_1 $.
$\Cox$
\end{pr}

\noindent{Proof} of Theorem~\ref{bdry nonconvergence}:  By hypothesis, 
condition (\ref{lnm}) holds,
and hence~(\ref{ineq}) holds for some choice of $\pp , k$ and $b$ which
are fixed hereafter.  Pick $\ee$ and $\nbd_1$ according to 
Proposition~\ref{1+b}.  Apply Proposition~\ref{convex} to $\nbd_1$
and $\pp$ with $\dd = 1\wedge (1+b/2)/(1+b/4)-1$ to obtain a neighborhood 
$\nbd$ of $\pp$ with the appropriate properties.  Temporarily fixing $n$, 
define the event $\B_n$ by $\VV (i) \in \nbd$ for all 
$n \leq i \leq (1+\dd) n$.  Define stopping times $\{ \tau_{i,r} \}$
and a family of Bernoulli random variables $\{ B_{i,r} \}$ as follows.  

\begin{quotation} \renewcommand{\baselinestretch}{1.0}\large\normalsize
\noindent{Let} $\tau_{i,r} \leq \infty$ be the $r^{th}$ time after 
$n$ that $Y_j = i$,
so formally $\tau_{i,0} = n$ and $\tau_{i,r+1} = \inf \{j > \tau_{i,r} :
Y_j = i\}$.  Let $B_{i,r} $ be independent and Bernoulli with 
\begin{equation} \label{bernoulli}
\P (B_{i,r} = 1 ) = \R_{ki} v_k (n) / (1+\dd ) N_i (n)
\end{equation}
and coupled to the variables $\{Y_i\}$ so
that if $B_{i,r} = 1 $ and $\tau_{i,r} \leq (1+\dd) n$ then $Y_{1+\tau_{i,r}}
= k.$   
\renewcommand{\baselinestretch}{1.3}\large\normalsize
\end{quotation}
To verify that this construction is possible, check that the probability
of a transition from vertex $i$ to vertex $k$ never drops below the quantity
in~(\ref{bernoulli}): 
\begin{eqnarray*}
\P (Y_{1+\tau_{i,r}} = k \| \F_{\tau_{i,r}} ) & \geq &  
(n / \tau_{i,r}) \R_{ki} v_k (n) / N_i (n) \\
& \geq & (1 / (1 + \dd ))\R_{ki} v_k (n) / N_i (n) 
\end{eqnarray*}
for $\tau_{i,r} < (1+\dd) n$.  

Now consider the subcollection $A \deq\{(i,r) : r \leq \dd n(p_i -\ee )\}$.
By Proposition~\ref{convex}, $\tau_{i,r} \leq (1+\dd) n$ whenever
the event $\B_n$ holds.  Meanwhile,
$$\E (\sum_{\alpha \in A} B_{\alpha} ) = \sum_i \dd n (p_i - \ee ) \R_{ki} 
   v_k (n) / (1+ \dd ) N_i (n).$$  
By Proposition~\ref{1+b}, this quantity is at least $\dd (1 + b/2)S_k (n) 
/ (1 + \dd )$ which is at least $\dd (1+b/4) S_k (n)$ by choice of $\dd$.  
Apply Proposition~\ref{pointmass} to the collection $\{ B_\alpha :
\alpha \in A \}$, with $b$ replaced by $b/4$ and $\ee_1$ to be chosen later 
to obtain a value for $L_0$.  Now calculate the conditional expectation
$\E (\ln (v_k ((1+\dd)n)) \| \F_n , S_k (n) > L_0)$.
By Proposition~\ref{pointmass} and the coupling,
\begin{eqnarray*}
&& \P (\B_n \mbox{ and } S_k ((1+\dd)n) - S_k (n) \geq \dd (1+b/8) S_k (n)  
   \| \F_n , S_k (n) > L_0) \\[2ex]
& > & \P (\B_n \| \F_n , S_k (n) > L_0) - \ee_1 .
\end{eqnarray*}
When $S_k ((1+\dd)n) - S_k (n) \geq \dd (1+b/8) S_k (n)$, it follows 
that $v_k ((1+\dd)n) \geq v_k (n) (1 + b/8(1+\dd)) \geq v_k (n) (1+b/16)$.
Therefore
\begin{eqnarray}
&& \E (\ln (v_k ((1+\dd)n)) \| \F_n , S_k (n) > L_0) \nonumber \\[2ex]
& \geq & (\P (\B_n \| \F_n , S_k (n) > L_0) -\ee_1) \ln ((1+b/16)v_k (n)) \\ 
&& +(1-\P (\B_n \| \F_n , S_k (n) > L_0) + \ee_1) \ln(v_k (n) /(1+\dd))
  \nonumber  \\[2ex] 
& \geq & \ln (v_k (n)) + \ln (1+b/32) - K \P (\B_n^c \| \F_n , S_k (n) 
   > L_0) \label{qwe}
\end{eqnarray}
for $K \ln (1+\dd) (1+b/16)$, when $\ee_1$ is sufficiently
small.  To conclude from this that $\P (\VV (n) \rightarrow \pp 
\mbox{ and } S_k (n) > L_0 \mbox{ for some } n) = 0$, write $T(n) =
(1+\dd)^n L_0$, $\G_n = \F_{T(n)}$, $X_n = \ln (v_k (T(n)))$, 
$\beta = c/2K$, $T = \inf \{ n : \VV (i)
\notin \nbd \mbox{ for some } L_0 \leq i \leq T(n) \}$ and calculate
\begin{eqnarray*}
\E X_{n \wedge T} & = & X_0 + \sum_{i=0}^{n-1} \E (1_{T > i} 
   (X_{i+1} - X_i) \| \G_i) \\[2ex]
& \geq & X_0 + \sum_{i=0}^{n-1} \E 1_{T > i} \left [ (c-K\beta) \cdot
    1_{\P (T = i+1 \| \G_n) \leq \beta} - K \cdot
    1_{\P (T = i+1 \| \G_n) > \beta} \right ] \\
&& \mbox{ by equation~\ref{qwe}} \\[2ex]
& \geq & X_0 + \sum_{i=0}^{n-1} \E 1_{T > i} \left [ (c-K\beta) \cdot
    (1 - \beta^{-1} \P (T = i+1 \| \G_n , T > i)) \right. \\
&&  \left. - K \cdot
    \beta^{-1} \P (T = i+1 \| \G_n , T > i) \right ] \\[2ex]
& \geq & X_0 + \sum_{i=0}^{n-1} (c-K\beta) 1_{T > i} - \beta^{-1} 
   (c+K-K\beta) \P (T = i+1) \\[2ex]
& \geq & X_0 + n (c-K\beta) \P (T > n) - \beta^{-1} 
   (c+K-K\beta) . 
\end{eqnarray*}
Since $c-K\beta$ was chosen to be positive, $\P (T > n)$ must go
to zero, showing that $\VV (n) \rightarrow \pp$ and $S_k (n) > L_0$
eventually is impossible.  

Finally, to show that $\P (\VV (n) \rightarrow \pp 
\mbox{ and } S_k (n) \leq L_0 \mbox{ for all } n) = 0$, note that
since $N_k (\pp ) > 0$, there is a sufficiently small neighborhood
$\nbd$ of $\pp$ for which $\P (Y_{i+1} = k \| \F_i , \VV (n) \in \nbd )$ 
is always at least a constant times $n^{-1}$.  Borel-Cantelli
implies that $k$ is visited infinitely often whenever
$\VV (n)$ remains in $\nbd$, and this finishes the proof of 
Theorem~\ref{bdry nonconvergence}.    $\Cox$ 

\section{Proof of Theorem~{\protect \ref{interior nonconvergence}}} 
\label{pf nonconverge2}

Begin with a proof of the equivalences:
\begin{eqnarray*}
&& \pp  \mbox{ fails to be a maximum for } H \\
& \Leftrightarrow & \R \mbox{ has more than one positive eigenvalue }\\
& \Leftrightarrow & \DT \mbox{ has a positive eigenvalue } 
\end{eqnarray*}

\noindent{The} matrix $\R$
can be viewed as a symmetric bilinear form whose quadratic form gives $H$ 
when restricted to $\simp$.  Let $\TD = \simp - \pp$ be the
translation of $\simp$ containing the origin.  For $\ww \in \TD$,
$$\R (\ww ,\pp ) = \ww^T \R \pp = \ww \cdot \lambda \cdot (1, \ldots , 1) = 0$$
where $\lambda$ is the common value of the $N_i$.  Then
\begin{equation} \label{W}
\R (\ww + c\pp , \ww + c\pp ) = \R (\ww ,\ww ) + \R (c\pp , c\pp )= 
\R|_{\TD}(\ww ) + c^2 \lambda
\end{equation}
so the quadratic form $\R (\vv , \vv )$ decomposes into the sum of $\R|_\TD$ 
and a positive form on the one-dimensional subspace spanned by $\pp$.  Then
$\R$ has precisely one more positive eigenvalue than the quadratic form
$R|_\TD$.  But equation~(\ref{W}) with $\ww = \vv - \pp$ shows that $H(\vv )
= \R|_\TD (\vv -\pp ) + \lambda$ so $H$ has a strict maximum at $\pp$ if and
only if $\R|_\TD$ has a strict maximum at the origin.  Since $\R$ has no 
zero eigenvalues, $R|_\TD$ will have a strict maximum when it has a 
maximum, which happens when it has no positive eigenvalues.    

For the second equivalence, note that $\ppi$ is smooth on the interior 
of $\simp$, so $\DT$ exists.  Let $\TTT$ be the operator on $\Real^d$
whose matrix in the standard basis is given by
$$\TTT_{ij} = \R_{ij} p_i / \lambda .  $$
I claim that $\TTT = \DT$ on $W$.  Indeed, using Definition~(\ref{eqpi})
to define $\ppi$ on all of $\Real^d$ and differentiating shows that
the matrix representation for $\DT$ is given by
\begin{eqnarray*}
[\DT]_{ij} & = & \left. {\partial \over \partial e_j} (\pi (\vv ))_i \right 
  |_{\vv = \pp} - \delta_{ij}  \\
& = & \left.  {\partial \over \partial e_j} {v_i N_i \over H(\vv )} \right 
  |_{\vv = \pp} - \delta_{ij} \\
& = & \left. \left ( {\R_{i,j} v_i \over H} + {\delta_{ij} N_i \over H} - 
  {v_i N_i \partial H /\partial e_j \over H^2} \right ) \right |_{\pp} - 
  \delta_{ij} \\
& = & \R_{ij} p_i / \lambda - 2 p_i 
\end{eqnarray*}
(using the fact that all the $N_i$ have a common value $\lambda = H(\pp )$ and
the identity ${\displaystyle {\partial H \over \partial e_j} = 2N_j}$).  
Then the matrices for $\TTT$ and $\DT$ differ by a matrix with constant
rows, hence define the same operator on $W$.  Now let
$diag(\pp )$ be the diagonal matrix with $i,i$ entry equal to $p_i$ and 
observe that $\TTT = diag(\pp ) \R / \lambda $.  
Since $\R$ is symmetric and $diag(\pp )$ is positive definite, 
$\TTT$ must be diagonalizable with real eigenvalues and has the same 
signature as $\R$ (see \cite[Theorem 6.23 and 6.24 page 232]{Or}).  Since 
$\TTT$ has $\pp$ as a positive eigenvalue and $\TD$ as an 
invariant subspace, it has one more positive eigenvalue than $\DT$ and the 
conclusion follows.   $\Cox$

To finish proving Theorem~\ref{interior nonconvergence}, 
it remains to show that $\vv (n)$ cannot converge to an interior point 
where $\DT$ has a positive eigenvalue.  The method of
proof is from \cite{Pe3}, the first step being a 
construction of a scalar function which measures ``distance from 
$\pp$ in an unstable direction'' (\cite[Proposition 3]{Pe3}).
\begin{lem} \label{quoteeta}
Under the assumptions of Theorem~\ref{interior nonconvergence},
suppose that $\DT$ has a positive eigenvalue.  Then
%(
there is a function $\eta$ from a neighborhood of $\pp$ to $[0,\infty )$
such that $\D_\vv \eta (\ppi (\vv ) - \vv) \geq k_1 \eta (\vv )$
in a neighborhood of $\pp$ for a constant $k_1 > 0$.  Furthermore,
$\eta$ is the square root of a smooth function (whose gradient necessarily
vanishes whenever the function vanishes) but whose second partials
are not all zero (thus $\eta$ is not differentiable where it vanishes). 
It follows from this that $\eta$ is Lipschitz and that $\eta (\vv + \ww ) 
\geq \eta (\vv ) + \D_\vv \eta (\ww ) + k_2|\ww |^2$ in a neighborhood
of $\pp$, where $\D_\vv \eta (\FF )$
may be any of the support hyperplanes to the graph of $\eta$
at points where $\eta$ vanishes.  $\Cox$
\end{lem} 

Use this lemma and a sequence of appropriately chosen stopping
times to convert questions about convergence of $\VV (n)$
into questions about the convergence of a scalar stochastic process.
To do this, fix a neighborhood $\nbd$ of $\pp$ in which all 
coordinates are bounded away from zero.  Let $L(\vv)$ be the
mean recurrence time to state 1 for the Markov chain $M(\vv)$
and let $L_{\max} = \sup_{\vv \in \nbd} L(\vv)$.  Pick $N_0
> 2L_{\max}$ and define $\sigma_0 = \inf \{k \geq N_0
: Y_k = 1 \}$ and $\sigma_{n+1} = \inf \{ k > \sigma_n : Y_k = 1\}$ 
to be the successive hitting times for state 1.  
Let $\tau = \inf \{ k \geq N_0 : \VV (k) \notin \nbd \}$ and let
$\tau_i = \tau \wedge \sigma_i$.  For the remainder of the section,
let $\Eb$ and $\Pb$ denote conditional expectation and conditional
probability with respect to $\F_{\tau_n}$.  The following 
facts are elementary. 

\begin{pr}  \label {tails} ~ \\
\begin{quotation}
$(i)$ The distribution of $\tau_{n+1} - \tau_n$ has finite conditional
expectation and variance.  Specifically,
$$\Pb (\tau_{n+1} - \tau_n \geq k + 1) < e^{- \alpha k}$$ 
for some $\alpha > 0$.

$(ii)$ For any $\ee > 0, N_0 , $ there is a constant $c_1$ such that 
$$ \P (n \leq \tau_n \leq c_1 n \mbox{ for all } 
  \tau_n \leq \tau ) > 1 - \ee . $$

$(iii)$ For any $\ee >0, \gamma < 1$, $N_0$ and $r$ may be chosen large enough 
so that
$$ \P (\tau_{n+1} - \tau_n \geq n^{1 - \gamma} \mbox{ for some } n \geq r) <
     \ee . $$
\end{quotation}
$\Cox$
\end{pr}

Let $\UU (n) = \VV (\tau_n )$, let $S_n = \eta (\UU (n))$ and let
$X_n = S_n - S_{n-1}$.  The following estimate shows that the 
expected increment in $\UU$ from time $n$ to $n+1$ is close to
the value given by the Markov approximation.

\begin{pr} \label{stationary}
For any $n > 0$,
\begin{equation} \label{old5.21}
\left | \Eb (\UU (n+1) - \UU (n)) - L(\UU (n)){\pi (\UU (n)) - 
   \UU (n) \over \tau_n} \right | = O ( {\tau_n}^{-2})  . 
\end{equation}
\end{pr}

\noindent{Proof:}  Couple the process $\{ Y_i : i \geq \tau_n \}$ to a 
Markov chain $Y_i'$ with $Y_{\tau_n} = 1$ and transition matrix  
$M(\UU (n))$ in such a way that the two processes remain identical
for as long as possible.  Define $\VV' , \SS', \tau'$ and $\UU'$ analogously
to the unprimed variables.  Establish first that
\begin{equation} \label {zeqc}
\Eb | \UU (n+1) - \UU' (n+1)| = O ( {\tau_n}^{-2} ) .
\end{equation}
To see this, observe that since transition probabilities for $Y$ and $Y'$
differ by at most $k / \tau_n$ at time $\tau_n + k$, the conditional 
probability of the two processes uncoupling before time $\tau_{n+1}$ is 
at most
\begin{equation} \label{za}
\sum_{k \geq 0} \Pb (\tau_{n+1} - \tau_n > k ) k/ \tau_n \leq e^{-\alpha}
  / (1-e^{-\alpha})^2 \tau_n
\end{equation}
according to Proposition~\ref{tails}~$(i)$.  On the other hand, 
$\E ( \tau_{n+1}' - (\tau_n + k) \| \F_{k + \tau_n})$ and
$\E ( \tau_{n+1}' - (\tau_n + k) \| \F_{k + \tau_n})$ are bounded
by $L_{\max}$ and $(1-e^{-\alpha})^{-1}$ respectively on the event
of the uncoupling occurring at time $k + \tau_n$, which implies that
\begin{equation} \label{zeqb} \begin{array} {ll}
& \Eb (|\vv (\tau_{n+1})-\VT (\tilt )| \;\|\mbox{ uncoupling before }
  \tau_{n+1}) \\[2ex]
\leq & \sup_k \Eb (|\vv (\tau_{n+1}) - \vv (\tau_n )| + |\VT (\tilt ) - 
  \vv (\tau_n )| \\[2ex]
& | \mbox{ uncoupling occurs at }\tau_n +k) \\[2ex]
  \leq & \sup_k (1 / \tau_n ) (\Eb (\tau_{n+1} - \tau_n -k+1 \\[2ex]
& + \tilt - \tau_n -k+1) \;| \mbox{ uncoupling occurs at }\tau_n +k) \\[2ex]
  \leq & (1 / \tau_n ) (L_{\max} + 1/(1-e^{-\alpha}) +2 ).
\end{array} \end{equation}
Combining~(\ref{za}) and~(\ref{zeqb}) gives~(\ref{zeqc}).

The quantity $\UU (n+1) - \UU (n)$ in the LHS of equation~(\ref{old5.21})
may now be replaced by the quantity $\UU' (n+1) - \UU (n)$, since
the two are within $O(\tau_n^{-2})$ in expectation.  Since $Y'$ is
a Markov chain, the following identity holds:
\begin{equation} \label{markovslln}
\Eb (\SS (\tau_{n_1} ) - \SS ( \tau_n ) ) = L(\UU (n)) \pi (\UU (n)).
\end{equation}
Component by component, we then have
\begin{eqnarray*}
& & \Eb ( U_i' (n+1) - U_i' (n)) \\[2ex]
& = & \Eb ( S_i' (\tau_{n+1}') / \tau_{n+1}' - S_i (\tau_n ) /\tau_n ) \\[2ex]
& = & \Eb \left ( {1 \over \tau_n} (S_i' (\tau_{n+1}) - S_i ( \tau_n ) - 
   (\tau_{n+1}' - \tau_n ) S_i (\tau_n ) / \tau_n ) \right ) - Q \\ 
& = & {1 \over \tau_n} L(\UU (n)) [\pi_i (\UU (n)) -U_i (n)] - Q 
\end{eqnarray*}
according to (\ref{markovslln}), where 
$$Q = \Eb \left ( {\tau_{n+1}' - \tau_n \over \tau_{n+1}' } \; \cdot \; 
   { S_i' ( \tau_{n+1}) - S_i (\tau_n) - (\tau_{n+1}' - \tau_n)S_i (n) 
   / \tau_n \over \tau_n} \right ) .$$
The denominator of $Q$ is ar least $\tau_n^2$ and the numerator is bounded
by the product of two geometric random variables according to
Proposition~\ref{tails}, so $|Q| = O(\tau_n^{-2})$ and the proposition
is proved.    $\Cox$

Use this estimate to prove the following proposition, which 
together with the Lemma~\ref{method} proves 
Theorem~\ref{interior nonconvergence}.

\begin{pr} \label{last}
Let $S_n$ and $X_n$ be defined from $\VV$ as above.  Let
$\nbd$ remain fixed as in the paragraph before 
Proposition~\ref{tails}.  For any $\ee > 0$ there are constants
$b_1 , b_2 , c > 0$ and $\gamma > 1/2$ and an
$N$ such that whenever $N_0 > N$ then 
\begin{equation} \label{one criterion}
\P (\B \| \F_{N_0}) > 1-\ee ,
\end{equation}
where $\B$ is the event that either 
equations~(\ref{lem1})~-~(\ref{lem4}) are satisfied for all $n > N_0$
or else $\VV (n)$ at some point leaves $\nbd$.
\begin{equation} \label{lem1}
\Eb ({X_{n+1}}^2 + 2 X_{n+1} S_n ) \geq b_1 / n^2
\end{equation}
\begin{equation} \label{lem2}
\Eb (X_{n+1} S_n  \one_{S_n > c/n} ) \geq 0
\end{equation}
\begin{equation} \label{lem3}
\Pb (|X_{n+1}| \leq 1 / (n+1)^{\gamma}) = 1
\end{equation} 
\begin{equation} \label{lem4}
\Eb ({X_{n+1}}^2 ) \leq b_2 / n^2
\end{equation}
\end{pr}

\begin{lem} \label{method}
If~(\ref{one criterion}) holds for a nonnegative stochastic process 
$S_n = S_0 + \sum_{i=1}^n X_i$, then $\P (S_n \rightarrow 0) = 0$.
\end{lem}

Lemma~\ref{method} is a variant on an argument from \cite{Pe3},
whose proof can be outlined as follows.

First assume that~(\ref{one criterion}) holds with $\ee = 0$, i.e. that
(\ref{lem1}) - (\ref{lem4}) hold almost surely,
and show in the following three
steps (A)-(C) that $\P (S_n \rightarrow 0) = 0$.  Let $k$
be any positive real number less than $\sqrt{b_1 /2}$ and
without loss of generality restrict $n$ to be at least $4c^2 /k$ so that
$k / 2\sqrt{n} > c/n$ and $c$ is the constant in condition~(\ref{lem2}).

\vspace{.2in}
\noindent{(A)} Claim: given any $S_n$, the probability of finding $S_M > 
k / \sqrt{n}$ for some $M \geq n$ is at least $1/2$. 
\begin{quotation}
Proof: Assume without loss of generality that $S_n < k/\sqrt{n}$.  Let 
$\sigma$ be the first $i \geq n$ for which $S_i > k/\sqrt{n}$.  
Then for any $M>n$,
\begin{eqnarray*}
&& \Eb (S_{\sigma \wedge M}^2) \\[2ex]
& = & S_n^2 + \sum_{i = n}^{M-1} \Eb (S_{\sigma \wedge (i+1)}^2 - 
      S_{\sigma \wedge i}^2 ) \\[2ex]
& = & \S_n^2 + \sum_{i = n}^{M-1} \Eb (\one_{\sigma > i} (X_{i+1}^2 +
      2 X_{i+1} S_i )) \\[2ex]
& \geq & \Pb (\sigma > M) \sum_{i = n}^{M-1} b_1 / i^2 \\
&& \mbox{ by } (\ref{lem2}) \\[2ex]
& \geq & \Pb (\sigma > M) b_1 / n .
\end{eqnarray*}
But condition~(\ref{lem3}) implies that $S_{\sigma \wedge i}$ never
gets much more than $k/\sqrt{n}$ and since $k^2 <
2 b_1$ this forces $\Pb (\sigma > M) < 1/2$ and the claim is proved. 
\end{quotation}
\vspace{.2in}
\noindent{(B)} Claim: given that $S_n  > k / \sqrt{n}$ the probability that 
$S_M$ will never return to the interval $x <k/\sqrt{n}$ for $M>n$ is at least
$a = 4b_2 / (4b_2 + k^2)$ .
\begin{quotation}
Proof: Assume $S_n > k/\sqrt{n}$.  Let $\sigma$ be the first $i>n$ for which
$S_i < k / 2\sqrt{n}$.  By condition~(\ref{lem2}) and the fact that
$S_{(\sigma - 1) \wedge i} > k / 2\sqrt{n} > c/n$, the sequence 
$S_{\sigma \wedge i}$ is a submartingale.  Decompose this into a mean-zero 
martingale and an increasing process.  Summing equation~(\ref{lem4}) shows 
the variance of the martingale to be bounded in $L^2$ by $b_2 / n$.  
Then by using the
one-sided Tschebysheff estimate $\P (f- \E f < -s) \leq \var (f) / 
(\var ( f) + s^2 )$, the probability that the martingale ever reaches
the interval $[-\infty , -k / 2 \sqrt{n})$ is at most 
$4 b_2 / (4 b_2 + k^2)$.  The martingale is a lower bound for the
submartingale so the claim is proved.
\end{quotation}
\vspace{.2in}
\noindent{(C)} If $S_n$ converges to 0 with non-zero probability, then there is an $n$
which can be chosen arbitrarily large and an
event $\A \in \Fn$ for which $\P (S_n \rightarrow 0 \| \A)$ is arbitrarily
close to 1.  When it is greater than $1 - a/2$, this contradicts (A) and (B).

\vspace{.2in}
Now assume~(\ref{one criterion}) 
instead of~(\ref{lem1}) - (\ref{lem4}).
For any $N_0$, let $\sigma$ be the first $n \geq N_0$ for which $\VV (n)$ exits 
$\nbd$ or one of the conditions~(\ref{lem1}) -~(\ref{lem4}) is violated;  
$\sigma$ is a stopping time since the conditions are $\Fn$-measurable.
Let $\{{X^*}_n , {S^*}_n \; : \; n \geq N_0 \}$ 
be any process that always satisfies~(\ref{lem1})~-
(\ref{lem4}) and is coupled to the process $\{X_n , S_n \; : \; n>N_0 \}$ 
so that the two processes are equal for $n \leq \sigma$.
Since $S_n^*$ cannot converge to $\pp$, $S_n \rightarrow \pp$ implies
$\sigma < \infty$.  For $\ee > 0$ let $N_0$ be chosen as 
in~(\ref{one criterion}).
Then with probability at least $1 - \ee$ either $S_n$ does not converge
to $\pp$ or $\VV (n)$ exits $\nbd$.  Thus the probability of $\VV (n)$
converging to $\pp$ without ever exiting $\nbd$ after time $N_0$ is at
most $\ee$.  Since $\ee$ is arbitrary, it follows that $\P (\VV (n)
\rightarrow \pp) = 0$.    $\Cox$

The last step in the proof of Theorem~\ref{interior nonconvergence}
is to establish Proposition~\ref{last}.
For any $\ee > 0$ and $\gamma < 1$, condition~(\ref{lem3}) may be satisfied
by choosing $N_0$ at least as large as the $N_0$ 
in Proposition~\ref{tails} $(iii)$ (using the fact that $\eta$ is Lipschitz).
Also,~(\ref{lem4}) follows
directly from Proposition~\ref{tails} $(i)$ for any $\ee$.
To prove~(\ref{lem2}), let $\ee > 0$ and use the bounds on $\tau_n$
from Proposition~\ref{tails} $(ii)$ to get
\begin{eqnarray}
\Eb (X_{n+1}) & = & \Eb (S_{n+1}) - S_n \nonumber \\
& = & \Eb ( \eta (\UU(n) + [\UU (N+1) - \UU (n)] )) - S_n \nonumber \\
& \geq & \Eb ( \eta ( \UU(n)) + \D_{\UU(n)} \eta 
  [\UU (n+1) - \UU (n)] + O|\UU (n+1) - \UU (n)|^2 ) - S_n \nonumber \\
&& \mbox{ by Lemma~\ref{quoteeta}} \nonumber \\
& = & \D_{\UU (n)} \eta \Eb [\UU (n+1) - \UU (n)] + \Eb (O|\UU (n+1) - 
    \UU (n)|^2 ) ) \nonumber \\
& = & \D_{\UU (n)} \eta \Eb \left [ {L(\UU(n)) \over \tau_n }
   (\pi - I ) \UU (n) + O(\tau_n^{-2}) \right ] + \Eb ( O|\UU (n)-\UU(n)|^2 )) 
    \nonumber \\
&& \mbox{ by Proposition~\ref{stationary}} \nonumber \\
& = & {L(\UU (n)) \over \tau_n } \D_{\UU (n)} \eta ((\pi - I)
    (\UU (n))) + O(\tau_n^{-2}) \nonumber \\
&& \mbox{ since $\UU (n+1) - \UU (n)$ is of order $\tau_n^{-1}$ 
    and $\eta$ is Lipschitz} \nonumber \\
& \geq & {k_1 L (\UU (n)) \over \tau_n} \eta (\UU (n)) + O(\tau_n^{-2})
    \nonumber \\
& \geq & {c_1 S_n \over n} - {c_2 \over n^2} \label{dominate}  \\
&& \mbox{ for some $c_1 , c_2 > 0$ by }
 \mbox{ Proposition}~\ref{tails} ~(ii)~ \mbox{ with probability }
  1-\ee . \nonumber 
\end{eqnarray}
Thus there is a constant $c = c_2 / c_1$ such
that for $S_n > c / n$ the first term of~(\ref{dominate}) dominates.  
Hence~(\ref{lem2}) is true with probability at least $1 - \ee$.

Finally, to show (\ref{lem1}), note that it suffices to show that
$\Eb (X_{n+1}^2 ) \geq c_3 / n^2$ for some $c_3$, assuming $\tau_n
\leq \tau$.  For, in the case that $S_n > c_2 / c_1 n$,
(\ref{lem2}) holds, implying~(\ref{lem1}), while if 
$S_n \leq c_2 / c_1$, (\ref{dominate}) is at least
$-2 c_2 \,n^{-2}$ and the second term on the left hand side
of~(\ref{lem1}) is at least $-4 c_2^2 / c_1 n^{-3}$ and for large
enough $n$ this is dwarfed by the $\Eb (X_{n+1}^2 )$ term.  

Now a moment's thought shows that $\Eb (X_{n+1}^2 )$ must be at least
order $n^{-2}$: from the nonvanishing second partials of $\eta^2$,
it follows that there is a unit vector $\ww \in W$ such that
$|\eta (\vv + r \ww) - \eta (\vv) > C r$  for some positive $C$
uniformly in $\vv$ in a neighborhood of $\pp$.  There exists a
positive multiple of $\ww$ and a fixed sequence
of sites $\{ 2 , \ldots , d \}$, such that if these are the
sites visited between times $\tau_n$ and $\tau_{n+1}$, then
$\UU (n+1) - \UU (n)$ will be arbitrarily close to this
multiple of $\ww$.  This sequence of visits happens with positive
probability, so~(\ref{lem1}) holds, establishing Proposition~\ref{last}
and Theorem~\ref{interior nonconvergence}.   $\Cox$

\section{Examples and further questions}

\noindent{Example 1}:  Suppose $\R_{ij} = 1 - \dd_{ij}$.  The critical 
set $\C$ contains just the centroids of the faces, and the
degeneracy set $\C_0$ is empty, so by Corollary~\ref{converge to point},
$\VV (n)$ converges almost surely to some point of $\C$.
It is easy to see that the centroids of all proper faces are linear
nonmaxima.  For example, if $\pp = (1/3 , 1/3 , 1/3,0,0,\ldots )$
then $N_i = 2/3$ for $i \leq 3$ and 1 for $i > 3$.  Thus 
Theorem~\ref{bdry nonconvergence} implies $\VV (n) \rightarrow
({1 \over d} , \ldots , {1 \over d})$ almost surely.

\noindent{Example 2}:  Here is an example where 
$\lim_{n \rightarrow \infty}$ is not deterministic.  Suppose 
$$ \R = \left [ \begin{array}{ccc} 3 & 1 & 1 \\ 1 & 2 & 4 \\ 1 & 4 & 2
   \end{array} \right ] \; . $$
All the minors of $\R$ are invertible and off-diagonal elements
nonzero, so Corollary~\ref{converge to point} applies and
$\VV (n)$ converges almost surely to a point of $\C$.  The interior
point $(1,1,1)\R^{-1} = (1/2 , 1/4 , 1/4)$ is unstable because
$\R$ has two positive eigenvalues, so the probability of convergence
there is zero.  The critical points in the middle of two of the
edges, $(1/3 , 2/3 , 0)$ and $(1/3 , 0, 2/3)$ are linear nonmaxima
as are the vertices $(0 , 1 , 0)$ and $(1 , 0 , 1)$, so the
probability of convergence to each of these points is zero by 
Theorem~\ref{bdry nonconvergence}.  On the other hand, $(1,0,0)$
is a local maximum for $H$ as is $(0,1/2,1/2)$, so by 
Theorem~\ref{th positive prob converge}, it follows that
$\P (\VV (n) \rightarrow (1,0,0)) = 1 - \P (\VV (n) \rightarrow
(0,1/2,1/2)) = a$ for some $0 < a < 1$.

\noindent{Example 3:}  Let $G$ be a finite abelian group and let
$T$ be a set of generators for $G$ closed under inverse.  Let
$\R$ be the incidence matrix for the Cayley graph of $(G,T)$.
By symmetry, the point $\pp = (1/|G| , \ldots , 1/|G|)$ is in $\C$.
The eigenvalues of $\R$ are just $\lambda (\chi) \deq \sum_{g \in T} 
\chi (g)$, as $\chi$ ranges over the characters of $G$.  If these
are all nonzero, then $\pp$ is the unique critical point
in the interior of $\simp$.  In this case, $\P (\VV (n) \rightarrow
\pp)$ is zero or not according to whether $\lambda (\chi) > 0$
for any nontrivial character $\chi$.  In fact it is easy to 
verify that $\P (\VV (n) \rightarrow \pp)$ is always zero 
or one when the principal minors of $\R$ are invertible, by 
checking that the negativity of $\lambda (\chi)$ for all 
nontrivial $\chi$ implies that each other critical point is 
a linear nonmaximum.

There are many natural unanswered questions about the behavior
of $\VV (n)$.  One could of course ask for rates of convergence, 
central limt behavior, etc., but I think it is more important
both from a mathematical and a modeling point of view
to try to extend the results already obtained so as to
cover all matrices $\R$.  For example, when $\R$ is a matrix
of all ones, every point of $\simp$ is critical so 
Theorem~\ref{th converge} says nothing, while comparison to 
a Polya urn model shows that $\VV (n)$
converges almost surely to a random point of $\simp$ with
an absolutely continuous distribution.  In general, when $\C$
has components larger than a point, one expects the motion
of $\VV$ inside a component to be martingale-like and hence
still converge to a single point, this time with a nonatomic
distribution.  Also, while
the symmetry assumption on $\R$ is vital to the proofs (since
it allows $\ppi (\vv )$ to be explicitly calculated) I do not
believe that it is actually necessary for the results.

\begin{conj} $\lim_{n \rightarrow \infty} \VV (n)$ exists almost
surely without any nondegeneracy assumptions on $\R$.  \end{conj}

\begin{conj} Theorem~\ref{th converge} holds whether or not $\R$ 
is symmetric.  Also, when $\R$ is not symmetric, there
is a function $H$ such that the first part of Lemma~\ref{Hincr}
holds and Theorem~\ref{interior nonconvergence} holds.  \end{conj}

\renewcommand{\baselinestretch}{1.0}\large\normalsize

\begin{flushleft}
\today
\end{flushleft}

\end{document}